\documentclass[11pt]{article}

\usepackage{amssymb,latexsym}
\usepackage{amsmath,amscd}
\usepackage{theorem}
\usepackage[all]{xy}

\title{\bf Generalized complex structures and Lie brackets}
\author{Marius Crainic \thanks{research supported by KNAW} \\[0.1cm]
         Department of Mathematics\\
        Utrecht University, P.O. Box 80.010, 3508 TA\\
        Utrecht, The Netherlands\\
        {\tt crainic@math.uu.nl}}
\date{}

\newcommand{\rmap}{\longrightarrow}

\newcommand{\SP} [1]     {{\left\langle {{#1}} \right\rangle}}

\newcommand{\pr}         {{\mathrm{pr}}}

\newcommand{\grd}         {\mathcal{G}}
\newcommand{\sour}        {\mathsf{s}}
\newcommand{\tar}         {{\mathsf{t}}}

\newtheorem{lemma} {Lemma} [section]
\newtheorem{proposition} [lemma] {Proposition}

\newtheorem{theorem} [lemma] {Theorem}

\newtheorem{definition}[lemma] {Definition}

\theorembodyfont{\rm} 

\newtheorem{remark}[lemma]{Remark}

\newenvironment{proof}{{\sc Proof:}}{{\hspace*{\fill} $\square$\\}}

\numberwithin{equation}{section}

\begin{document}

\maketitle

\begin{abstract}
We remark that the equations underlying the notion of generalized complex structure have
simple geometric meaning when passing to Lie algebroids/groupoids.
\end{abstract}

\tableofcontents

\section{Introduction}

Generalized complex structures have recently been introduced by N. Hitchin \cite{Hi} and further studied by M. Gualtieri \cite{Gu},
as a bridge between symplectic and complex geometry. When making the conceptual definition
more explicit, one arises at puzzling equations (see Section \ref{Generalized complex structures}). 
In this paper we remark that such complicated
explicit formulas have miraculously simple interpretations after allowing algebroids and groupoids into the
picture.
Since underlying a generalized complex structure there is a Poisson bivector, we use either as a guide, or 
explicitly, some basic principles from Poisson geometry. The first one we have in mind is that, when passing to the global objects behind Poisson structures,
one arises at symplectic structures (and similarly after restricting to leaves). Here, instead of ``symplectic'', please read
``{\it non-degenerate} Poisson''. Using the groupoid associated to the underlying Poisson structure, we show that the puzzling equations 
mentioned above translate into simple structures on the groupoid, giving rise to what we call here ``Hitchin groupoids''. In particular, 
we will see that, at the level of the groupoid, one has an induced {\it non-degenerate} generalized complex structure (but, in comparison with Poisson geometry, the situation is even better since the structure is more then non-degenerate). This is done in the last section, which starts with a short presentation of algebroids/groupoids, as well as of the main ingredient that
we use, namely the multiplicative 2-forms discussed in \cite{BCWZ} (however, the presentation here is self-contained). 
The first section is ``groupoid-free''; there we look at the definition of generalized complex manifolds (in terms of the
tensors involved), we introduce the generalized holomorphic maps (pointing out ``reduction'' in this context), and we have a closer look at 
generalized complex structures which are non-degenerate. We end the paper with a few remarks inspired from Poisson and Dirac geometry.

Our results can be viewed as ``integration results'', i.e. as  describing the global structure behind generalized complex structures.
Since underlying our discussion is the symplectic groupoid of the Poisson manifold, or the phase-space of the Poisson sigma model \cite{CaFe}, 
it would be interesting to find the relationship between our work and the 
``Hitchin sigma model'' 
\cite{final, Zu}. Another 
interesting (and probably related) question is to find the relationship between our integration results and the integration of complex Lie algebroids
proposed by A. Weinstein \cite{We}. Finally,  
this work is similar to \cite{new}, in the sense that we study generalized complex structures from the point of view of Poisson 
(but also Dirac) geometry.

\section{Generalized complex structures}
\label{section2}

\subsection{Generalized complex structures}
\label{Generalized complex structures}

Let $M$ be a (real) manifold, and let
\[ \mathcal{T}M:= TM\oplus T^*M ,\]
which we will call the generalized tangent bundle of $M$. It 
relates to the usual tangent bundle via the obvious projection, denoted:
\[ \pi: \mathcal{T}M\rmap TM ,\]
The space of sections of $\mathcal{T}M$ is endowed with a bracket, known as \textbf{the Courant bracket} \cite{Cou}, which extends the
Lie bracket on vector fields. Explicitly,
\[ [(X, \xi), (Y, \eta)]= ([X, Y], L_X\eta- L_Y\xi- \frac{1}{2}d(i_X\eta- i_Y\xi)).\]
Note that the Courant bracket is antisymmetric but it does not satisfy the Jacobi identity.

The generalized tangent bundle $\mathcal{T}M$ carries various 
other interesting structures, such as the canonical (fiberwise) symplectic form, or the
non-degenerate symmetric bilinear form $\langle \cdot, \cdot\rangle$:
\[ \langle (X, \xi), (Y, \eta) \rangle= \xi(Y)+ \eta(X) .\]
Although these will not play much role in this paper, they should be viewed as part 
of the picture. For more details on the relation between $\langle \cdot, \cdot\rangle$
and the Courant bracket, we refer to \cite{Cou, ca}.

By replacing the usual tangent bundle with $\mathcal{T}M$ and the Lie bracket with the Courant bracket,
many aspects of tensorial geometry can be carried out in this ``generalized'' setting. 
The famous Newlander-Nirenberg theorem enables us to apply this principle to complex structures. Let us first
recall that, for a bundle map $a:TM\rmap TM$ covering the identity, its torsion $\mathcal{N}_a$ is the $(1, 2)$-tensor:
\[ \mathcal{N}_a(X, Y)= [aX, aY]+ a^2[X, Y]- a([aX, Y]+ [X, aY]),\] 
and one says that $a$ is torsion-free if $\mathcal{N}_a=0$. With this,  the Newlander-Nirenberg theorem 
asserts that there is a 1-1 correspondence between complex structures on $M$ and almost complex structures $J$ 
on $M$ (i.e. vector bundle maps $J:TM\rmap TM$
satisfying $J^2= -Id$) which are torsion-free. Hence one arises \cite{Gu, Hi} at

\begin{definition} A \textbf{generalized complex structure on $M$} (g.c. structure on short) 
is a complex structure $\mathcal{J}$ on the bundle $\mathcal{T}M$
(i.e. a bundle map $\mathcal{J}: \mathcal{T}M\rmap \mathcal{T}M$ with $\mathcal{J}^2= -Id$), satisfying:
\begin{equation}
\label{J-torsion} 
[\mathcal{J}\alpha, \mathcal{J}\beta]- [\alpha, \beta]- \mathcal{J}([\mathcal{J}\alpha, \beta]+ [\alpha, \mathcal{J}\beta]) = 0,
\end{equation}
for all sections $\alpha, \beta$ of $\mathcal{T}M$. 
\end{definition}

To make the definition of generalized complex structures more explicit, we need to fix some notations:
for a two-form $\sigma$ on $M$ we denote by $\sigma_{\sharp}: TM\rmap T^*M$ the map $X\mapsto i_{X}\sigma$,
while for a bivector $\pi$ on $M$ we denote by $\pi^{\sharp}: T^*M\rmap TM$ the contraction with $\pi$, i.e.
defined by
\[ \beta(\pi^{\sharp}(\alpha))= \pi(\alpha, \beta) .\]
Also, we denote by $[\cdot, \cdot]_{\pi}$ the bracket on the space of 1-forms on $M$ defined by
(see also the last section):
\begin{equation}
\label{br-pi} 
[\xi, \eta]_{\pi}= L_{\pi^{\sharp}\xi}\eta- L_{\pi^{\sharp}\eta}\xi - d\pi(\xi, \eta) .
\end{equation}

Next, after spelling out (\ref{J-torsion}), reducing the (apparently longer) list of equations to a minimal one, and then
writing and organizing them in a convenient way, we arrive at the following (in local coordinates, the equations bellow also appear
in \cite{final}).


\begin{proposition}\label{first} A generalized complex structure $\mathcal{J}$ on $M$ is necessarily of type 
\begin{equation}
\label{J} 
\mathcal{J}= \left( \begin{array}{ll}
                a & \pi^{\sharp} \\
                \sigma_{\sharp} & - a^*  
        \end{array}
   \right) 
\end{equation}
where $\pi$ is a bivector on $M$, $\sigma$ is a two-form on $M$, and $a: TM\rmap TM$ is a bundle map.
Moreover, $\mathcal{J}$ is a generalized complex structure if and only if the following conditions are satisfied
\begin{enumerate}
\item[(C1)] $\pi$ satisfies the equation:
\[ \pi^{\sharp}([\xi, \eta]_{\pi})= [\pi^{\sharp}(\xi), \pi^{\sharp}(\eta)].\]
\item[(C2)] $\pi$ and $a$ are related by the following two formulas:
\begin{eqnarray}
 & a\pi^{\sharp} =  \pi^{\sharp}a^{*} \label{2.1}\\
 & a^*([\xi, \eta]_{\pi})  =  L_{\pi^{\sharp}\xi}(a^*\eta)- L_{\pi^{\sharp}\eta}(a^*\xi)- d\pi(a^*\xi, \eta)
\label{2.2}
\end{eqnarray}
\item[(C3)] $\pi$, $a$ and $\sigma$ are related by the following two formulas:
\begin{eqnarray}
 & a^2+ \pi^{\sharp}\sigma_{\sharp}= -Id \label{3.1}\\
 & \mathcal{N}_a(X, Y)= \pi^{\sharp}i_{X\wedge Y}(d\sigma) \label{3.2}
\end{eqnarray}
\item[(C4)] $\sigma$ and $a$ are related by the following two formulas:
\begin{eqnarray}
 & a^{*}\sigma_{\sharp} = \sigma_{\sharp}a\label{4.1}\\
 & d\sigma_a(X, Y, Z)= d\sigma(aX, Y, Z)+ d\sigma(X, aY, Z)+ d\sigma(X, Y, aZ) ,\label{4.2}
\end{eqnarray}
where $\sigma_a(X, Y)= \sigma(aX, Y)$.
\end{enumerate}
(the previous equations must hold for all 1-forms $\xi$ and $\eta$, and all vector fields $X$, $Y$ and $Z$).
\end{proposition}

Note that condition (C1) is equivalent to that fact that $\pi$ is a \textbf{Poisson bivector} (i.e. the bracket $\{f, g\}= \pi(df, dg)$
is a Lie algebra bracket on the $C^{\infty}(M)$). In particular, the distribution $\pi^{\sharp}(T^*M)\subset TM$ is integrable,
and its leaves (called the symplectic leaves) carry natural symplectic structures; also, 
the bracket $[\cdot, \cdot]$ restricts to a Lie bracket on $\mathfrak{g}_x= Ker(\pi^{\sharp}_{x})$, making it into a Lie algebra 
called the \textbf{isotropy Lie algebra} at $x\in M$. Of course, the first equations in (C2) and (C3) imply that $\mathfrak{g}_x$ 
is a complex Lie algebra (with the restriction of $a$ as the complex structure).

For latter use note also that if (\ref{J}) is a generalized complex structure, then so is
\begin{equation}
\label{opp-J}
\overline{\mathcal{J}}= \left( \begin{array}{ll}
                a & -\pi^{\sharp} \\
                -\sigma_{\sharp} & - a^*  
        \end{array}
   \right) ,
\end{equation}
which will be called \textbf{the opposite} of $\mathcal{J}$.

\begin{proof}
First of all, note that $\mathcal{J}$ must be orthogonal with respect to the inner product $\langle \cdot, \cdot\rangle$:
\[ \langle \mathcal{J}(\alpha), \mathcal{J}(\beta)\rangle = \langle \alpha, \beta\rangle \]
for all sections $\alpha$, $\beta$ of $\mathcal{T}M$. This follows from the $C^{\infty}(M)$-linearity  
of the equation (\ref{J-torsion}) with respect to $\beta$ and the formula:
\[ [\alpha, f\beta]= f[\alpha, \beta]+ L_{\pi(\alpha)}(f)\beta- \langle \alpha, \beta\rangle df,\]
for $f\in C^{\infty}(M)$. From this orthogonality property it immediately follows that $\mathcal{J}$ must be
of type (\ref{J}), while the condition $\mathcal{J}^{2}= -Id$ immediately implies that $\pi$, $a$ and $\sigma$ 
satisfy the linear equations appearing in (C2), (C3) and (C4). We claim that the remaining (differential) 
equations are equivalent to the integrability condition (\ref{J-torsion}). Let us look closer at (\ref{J-torsion}).
The equation has two components hence it gives us two equations. And this happens in each of the three cases
we have to consider. The first one is when $\alpha$ and $\beta$ are 1-forms, call them $\xi$ and $\eta$,
when the resulting equations are easily seen to be (C1) (for the first component) and the second one in (C2) (for the second component).
The next case is when $\alpha$ is a vector field $X$ and $\beta$ is a 1-form $\xi$. The first component of the resulting equation 
translates into:
\begin{equation}
\label{int1} 
[a(X), \pi^{\sharp}(\xi)]- a([X,  \pi^{\sharp}(\xi)])=  \pi^{\sharp}(L_{aX}(\xi)- L_{X}(a^{*}(\xi))).
\end{equation}
while the second component gives us:
\begin{eqnarray} 
 & \sigma_{\sharp}([\pi^{\sharp}(\xi), X])- L_{\pi^{\sharp}(\xi)}(\sigma_{\sharp}(X)) =  L_{X}(\xi)+ d \xi(\pi^{\sharp}\sigma_{\sharp}(X)) \nonumber  \\
 & + L_{a(X)}(a^{*}(\xi))- a^{*}(L_{a(X)}(\xi)- L_{X}(a^{*}(\xi))), \label{int2}
\end{eqnarray}
The third case is when $\alpha$ and $\beta$ are vector fields, call them $X$ and $Y$. The first component of the resulting equation is
\begin{eqnarray} 
 & [X, Y]- [a(X), a(Y)] + a([a(X), Y]+ [X, a(Y)]) =  \nonumber  \\
 & -\pi^{\sharp}(L_{X}(\sigma_{\sharp}(Y))- L_{Y}(\sigma_{\sharp}(X))+ d\sigma(X, Y)), \label{int3}
\end{eqnarray}
while the second component gives us:
\begin{eqnarray} 
 & -L_{a(X)}(\sigma_{\sharp}(Y))+ L_{a(Y)}(\sigma_{\sharp}(X))- d\sigma_{a}(X, Y) + \sigma_{\sharp}([a(X), Y]+ [X, a(Y)])= \nonumber  \\
 & a^{*}(L_{X}(\sigma_{\sharp}(Y))- L_{Y}(\sigma_{\sharp}(X))+ d\sigma(X, Y)). \label{int4}
\end{eqnarray}
Now, the equation (\ref{int1}) is easily seen to be equivalent to the the second equation in (C2) (again). Also, a straightforward 
computation shows that (\ref{int2}) and (\ref{int3}) are equivalent. Making use of the formula:
\begin{equation}
\label{differential} 
i_{X\wedge Y}(d\sigma)= L_{X}(i_Y(\sigma))- L_{Y}(i_X(\sigma))+ d(i_{X\wedge Y}\sigma)- i_{[X, Y]}\sigma,
\end{equation}
(\ref{int3}) is easily seen to be equivalent to the second equation of (C3). Finally, making
use of the Koszul-type formula
\begin{eqnarray} 
(d\sigma)(X, Y, Z) & =  L_{X}\sigma(Y, Z)+ L_{Y}\sigma(Z, X)+ L_{Z}\sigma(X, Y) \nonumber  \\
                   & -  \sigma([X, Y], Z)- \sigma([Z, X], Y)- \sigma([Y, Z], X), \label{Koszul}
\end{eqnarray}
a straightforward computation shows that (\ref{int4}) is equivalent to the second equation in (C4).
\end{proof}

Finally, let us point out the relation with Dirac geometry. 
Recall that a Dirac structure on $M$ is a sub-bundle $L\subset \mathcal{T}M$ of rank equal 
to the dimension of $M$, with the property that the canonical pairing $\langle \cdot, \cdot\rangle$ vanishes when restricted
to $L$, and which satisfies the involutivity condition $[\Gamma(L), \Gamma(L)]\subset \Gamma(L)$ (with respect to the Courant
bracket). The standard examples of Dirac structures are graphs of closed 2-forms or of Poisson bivectors on $M$. In our context,
condition (C2) is equivalent to saying that
\begin{equation}
\label{the-Dirac} 
L_{\pi, a}:= \{ (\pi^{\sharp}(\xi), a^{*}(\xi)): \xi\in T^*M \}
\end{equation}
is a Dirac structure on $M$. 

Complex Dirac structures are defined similarly, by replacing $\mathcal{T}M$ by its complexification
\[ \mathcal{T}_{\mathbb{C}}M= \mathcal{T}M\otimes \mathbb{C}= T_{\mathbb{C}}M\oplus T_{\mathbb{C}}^{*}M .\]

The main point of this discussion is that generalized complex structures $\mathcal{J}$ on $M$ can be identified with complex Dirac structures $L$ on $M$ with the property that $L\oplus \overline{L}= \mathcal{T}_{C}M$ \cite{Gu}. 
In one direction, given $\mathcal{J}$, $L$ is the $+i$-eigenspace of $\mathcal{J}$ on $\mathcal{T}M\otimes \mathbb{C}$. 
We see that, if $\mathcal{J}$ is represented
by (\ref{J}), then the associated $L$ is:
\[ L:= \{ (X-iX', \xi- i\xi'): X'= aX+ \pi^{\sharp}\xi, \xi'= \sigma_{\sharp}X- a^{*}\xi \}.\]

\subsection{The non-degenerate case}
\label{The non-degenerate case}

In this section we discuss the generalized complex structure which are non-degenerate in the sense that 
the associated bivector $\pi$ is non-degenerate (i.e. $\pi^{\sharp}$ is an isomorphism). Such structures 
arise after restricting to the symplectic leaves of $\pi$ or after passing to the global objects (see the 
next section). Also, they are particularly simple since some of 
the structure data (tensors and equations) are consequences of the non-degeneracy 
condition. To make this explicit, we first introduce some notations.

\begin{definition}\label{def-s+c} Given a 2-form $\omega$ on $M$, and a bundle map $a:TM\rmap TM$, we say that
$\omega$ and $a$ commute if $\omega_{\sharp}a= a^{*}\omega_{\sharp}$, i.e.:
\[ \omega(aX, Y)= \omega(X, aY),\]
for all $X, Y$-vector fields on $M$. In this case, we denote by $\omega_{a}$ the 2-form:
\[ \omega_a(X, Y)= \omega(aX, Y).\]
An \textbf{symplectic$+$complex structure} on $M$ consists of a pair $(\omega, J)$ with $\omega$-symplectic, and $J$-complex structure on $M$, which commute.
\end{definition}

For instance, in the (real) cotangent bundle of any complex manifold, the canonical symplectic form
together with the lifted complex structure is a symplectic$+$complex structure structure. 
Also, the co-adjoint orbits of a complex Lie algebra carry such structures. In our context, 
the symplectic$+$complex structures correspond to those 
non-degenerate generalized complex structures with $\sigma= 0$ (see below).

\begin{lemma}
\label{lemma-complex+sympl} 
Let $(M, \omega)$ be a symplectic manifold. Then 
\[ a\mapsto \omega_{a} \]
defines a 1-1 correspondence between $(1, 1)$-tensors $a: TM\rmap TM$ commuting with $\omega$ and 2-forms on $M$.
Moreover, $(\omega, a)$ is a symplectic$+$complex structure 
if and only if 
\[ d\omega_{a}= 0, \ \omega+ a^*\omega = 0.\]
\end{lemma}

\begin{definition}
A \textbf{Hitchin pair} on $M$ is a pair $(\omega, a)$ consisting of a symplectic form $\omega$
and a $(1, 1)$-tensor $a$ with the property that $\omega$ and $a$ commute and 
\[ d\omega_{a}= 0. \]
The twist of $(\omega, a)$ is the 2-form
\[ \sigma:= -(\omega+ a^*\omega) .\]
\end{definition}

Note that in this case $a$ is neither an almost complex structure nor torsion free. Instead, as it will follow from the computation below, one has:
\[ a^{2}= -Id- \omega_{\sharp}^{-1}\sigma_{\sharp},\ \mathcal{N}_a(X, Y)= \omega_{\sharp}^{-1}i_{X\wedge Y}(d\sigma) .\]


\begin{proposition}
\label{nondeg} 
There is a 1-1 correspondence between generalized complex structures 
\begin{equation} 
\mathcal{J}= \left( \begin{array}{ll}
                a & \pi \\
                \sigma_{\sharp} & - a^*  
        \end{array}
   \right) 
\end{equation}
with $\pi$ non-degenerate, and Hitchin pairs $(\omega, a)$. In this correspondence, $\pi$ is the inverse of $\omega$, and
$\sigma$ is the twist of the Hitchin pair $(\omega, a)$.
\end{proposition}

The remaining of this section is devoted to the proof of 
Lemma \ref{lemma-complex+sympl} and of Proposition \ref{nondeg}. We start with the proposition. 
At the level of the tensors (bivectors, 2-forms, $(1, 1)$-tensors), the correspondence is completely described in the statement, and 
we only have to take care of the equations they satisfy. 
The following is the well-known correspondence between non-degenerate Poisson bivectors and symplectic forms, which we include for completeness. 

\begin{lemma} If $\pi$ is a non-degenerate bivector on $M$, and $\omega$ is the inverse 2-form (defined by $\omega_{\sharp}= (\pi^{\sharp})^{-1})$,
then $\pi$ satisfies (C1) if and only if $\omega$ is closed.
\end{lemma}

\begin{proof} Apply (C1) to $\xi= i_{X}(\omega)$, $\eta= i_{Y}(\omega)$, where $X$ and $Y$ are arbitrary vector fields, then apply $\omega_{\sharp}$ to the resulting formula, and use (\ref{differential}).
\end{proof}

Next, we take care of (C2). 

\begin{lemma} Given a symplectic form $\omega$, $\pi$ the associated non-degenerate bivector (i.e. $\pi^{\sharp}= \omega_{\sharp}^{-1}$), 
and $a: TM\rmap TM$ a bundle map, then $\pi$ and $a$ satisfy (C2) if and only if $\omega$ and $a$ commute and $\omega_{a}$ is closed.
\end{lemma}

\begin{proof}
Clearly, (\ref{2.1}) is equivalent to $\omega$ and $a$ commuting. Next, to take care of (\ref{2.2}), 
we apply it to $\xi= i_{X}\omega$, $\eta= i_{Y}\omega$, where $X$ and $Y$ are arbitrary vector fields. Using
the $a^{*}(i_{X}(\omega))= i_{X}(\omega_a)$, and applying formula (\ref{differential}) to our closed $\omega$, we get:
\[ i_{[X, Y]}(\omega_{a})= L_{X}(i_{Y}(\omega_{a})- L_{Y}(i_{X}(\omega_a))+ d\omega_{a}(X, Y).\]
Hence, using again (\ref{differential}), we see that (\ref{2.2}) is equivalent to $\omega_a$ being closed.
\end{proof}

Next, the equation (\ref{3.1}), under the hypothesis that $\pi$ is non-degenerate, clearly forces
$\sigma=  - \omega- a^{*}\omega$,  where $\omega$ is the inverse of $\pi$. Hence we are left with proving that 
if $(\omega, a)$ is as in the statement of the proposition,
then $J$ given by (\ref{J}) with $\pi$ the inverse of $\omega$ and $\sigma=  - \omega- a^{*}\omega$
automatically satisfies equation (\ref{3.2}), (\ref{4.1}), (\ref{4.2}). The first one is taken care by the following general
computation:

\begin{lemma}
\label{torsion} 
For any 2-form $\omega$ and any $a:TM\rmap TM$ commuting with $\omega$, one has:
\[ i_{\mathcal{N}_a(X, Y)}(\omega)= i_{aX\wedge Y+ X\wedge aY}(d\omega_{a})- i_{aX\wedge aY}(d\omega)- i_{X\wedge Y}(d(a^{*}\omega)).\]
\end{lemma}

\begin{proof}
Using the Koszul-type formula (\ref{Koszul})
applied to each of the four terms appearing in the right hand side of the equation, and the fact that $\omega$ and $a$ commute,
we easily arrive to the left hand side. 
\end{proof}

Hence, under our hypothesis ($\omega$ and $\omega_{a}$-closed), the equation in the lemma implies that
\begin{equation}
\label{help1} 
i_{X\wedge Y}(d(a^{*}\omega))=  -i_{\mathcal{N}_a(X, Y)}(\omega).
\end{equation}
Hence $\sigma= -\omega- a^*\omega$ has:
\[ i_{X\wedge Y}(d\sigma)= i_{\mathcal{N}_a(X, Y)}(\omega),\]
and then, applying $\pi^{\sharp}$, we arrive at (\ref{3.2}). Next, equation (\ref{4.1}) is clear.
Finally, we look at (\ref{4.2}). Writing the equation as
\[ i_{X\wedge Y}(d\sigma_a)= i_{aX\wedge Y+ X\wedge aY}(d\sigma)+ a^{*}(i_{X\wedge Y}(d\sigma)),\]
since $\sigma= -\omega- a^*\omega$ and since $\omega$ and $\omega_{a}$ are closed, the desired equation is equivalent to:
\begin{equation}
\label{toprove} 
i_{X\wedge Y}d(a^*\omega_a)= i_{aX\wedge Y+ X\wedge aY}d(a^*\omega)+ a^*i_{X\wedge Y}d(a^*\omega).
\end{equation}
To compute the left hand side, we use Lemma \ref{torsion} applied to $\omega_a$ and we obtain:
\[ i_{X\wedge Y}d(a^*\omega_a)= i_{aX\wedge Y+ X\wedge aY}d(a^*\omega)- i_{\mathcal{N}_a(X, Y)}(\omega_a).\]
Using $i_{\mathcal{N}}(\omega_a)= a^*i_{\mathcal{N}}(\omega)$ and (\ref{help1}), we arrive at the desired equation 
(\ref{toprove}). This concludes the proof of Proposition \ref{nondeg}.

We now prove Lemma \ref{lemma-complex+sympl}.
The only part which needs some explanation is to show that if $a$ is almost complex then
$\mathcal{N}_a= 0$ if and only if $d\omega_{a}= 0$. 
If $\omega_{a}$ is closed, then Lemma \ref{torsion} together with the fact that $\omega$ is closed, $a^*\omega= - \omega$ and
$\omega$ is non-degenerate, implies that $\mathcal{N}_{a}= 0$. The converse follows again from a more general equation which is 
stated in the following lemma. 

\begin{lemma} If $J$ is a complex structure on $M$ then, for any 2-form $\omega$ on $M$ commuting with $J$ one has:
\begin{eqnarray*}
 & 2d\omega_{J}(X, Y, Z) = &  \\
 & d\omega(JX, Y, Z)+ d\omega(X, JY, Z)+ d\omega(X, Y, JZ)+ d\omega(JX, JY, JZ).
\end{eqnarray*}
\end{lemma}

\begin{proof} This follows from the general properties of complex forms, or can be 
checked directly, e.g. by using again the Koszul formula (\ref{Koszul}).
\end{proof}

\subsection{Generalized holomorphic maps}
\label{Generalized holomorphic maps and reduction}

Next, we discuss the notion of morphism between generalized complex manifolds. 
At first one may expect that, based on the interpretation of generalized complex structures
as complex Dirac structures (see the end of section \ref{Generalized complex structures}), one could use the
notion of morphism in Dirac geometry. However, there are two such notions: ``forward'' morphisms which are well-behaved
with respect to bivectors but not 2-forms,  and ``backward'' morphisms which are well behaved with respect to 2-forms
but not bivectors; see \cite{BCWZ, BO}. However, since generalized complex structures
are made up of both 2-forms and bivectors, both of these two notions would be quite restrictive. 
However, with the mind at the matrix representation (\ref{J}), there is a natural notion of morphism
which sits in between the two notions of forward and backward Dirac maps.

\begin{definition} Let $(M_i, \mathcal{J}_i)$, $i\in\{1, 2\}$, be two generalized complex manifolds, and let
$a_i, \pi_i, \sigma_i$ be the components
of $\mathcal{J}_i$ in the matrix representation (\ref{J}). A map $f: M_1\rmap M_2$ is called \textbf{generalized holomorphic} if 
$f$ maps $\pi_1$ into $\pi_2$, $f^*\sigma_2= \sigma_1$ and $(df)\circ a_1= a_2\circ (df)$.
\end{definition}

\begin{lemma}\label{fibers} If $f: M_1\rmap M_2$ is generalized holomorphic and $x\in M_2$ is regular for $f$,
then the fiber $\mathcal{C}= f^{-1}(x)$ is a complex manifold (with complex structure $a_1|_{\mathcal{C}}$). 
\end{lemma}

\begin{proof} Let $T\mathcal{C}= Ker(df)$ and since $(df)$ commutes with the $a_i$'s,
$a_{1}$ restricts to an $(1, 1)$-tensor on $\mathcal{C}$. Using the two equations in (C3) 
applied to $J_1$ and the fact that $f^*\sigma_2= \sigma_1$, we deduce that $a_{1}^2(X)= - X$
and $\mathcal{N}_{a_1}(X, Y)= 0$ if $X\in Ker(df)$. 
\end{proof}

With the mind at Poisson manifolds, it is very tempting to discuss reduction in our setting.

\begin{definition} A Hitchin realization of a generalized complex manifold $(M, \mathcal{J})$ consists
of a manifold $S$ endowed with a Hitchin pair $(\omega, a)$ and a smooth map $\mu: S\rmap M$ 
with the property that $\mu$ is generalized holomorphic. 
\end{definition}

Recall from the previous subsection that Hitchin pairs correspond to non-degenerate generalized complex structures.
A basic example of Hitchin realization is the inclusion of the symplectic leaves into $M$. 

We now discuss reduction. Let $\mu: S\rmap M$ be a Hitchin realization, and let $\pi, a, \sigma$ be the components of
$\mathcal{J}$. Assume that $x\in M$ is regular for $\mu$, and we look at the level set
\[ \mathcal{C}= \mu^{-1}(x) .\]
There is an induced infinitesimal action of the isotropy Lie algebra $\mathfrak{g}_x= Ker(\pi^{\sharp}_{x})$ on $\mathcal{C}$,
encoded into a Lie algebra map $\rho: \mathfrak{g}_x\rmap \mathcal{X}(\mathcal{C})$: for $\alpha\in \mathfrak{g}_x$,
$\rho(\alpha)$ is the unique solution of the moment map equation 
\[ i_{\rho(\alpha)}(\omega)= \mu^*(\alpha) .\]
We now assume that the quotient space $\mathcal{C}/\mathfrak{g}_x$ is smooth. The reduction from Poisson geometry 
says that $\mathcal{C}/\mathfrak{g}_x$ carries a canonical symplectic structure $\omega_{\mathcal{C}}$ (uniquely determined by
the condition that its pull-back to $\mathcal{C}$ coincides with the restriction of $\omega$ to $\mathcal{C}$).
The following shows that $\mathcal{C}/\mathfrak{g}_x$ carries a canonical symplectic$+$complex structure.

\begin{proposition} $\mathcal{C}/\mathfrak{g}_x$ carries a complex structure $J_{\mathcal{C}}$ unique with the
property that the projection $\mathcal{C}\rmap \mathcal{C}/\mathfrak{g}_x$ is a holomorphic map. Moreover, $J_{\mathcal{C}}$
commutes with the reduced symplectic structure $\omega_{\mathcal{C}}$.
\end{proposition}

\begin{proof} One can check directly that $J$ vanishes on the image of $\rho$, and then prove that $J$ can be pushed forward to
a $(1, 1)$-tensor on $\mathcal{C}/\mathfrak{g}_x$ with the desired properties. Alternatively, the same argument that shows that $\omega$ can be reduced also applies to $\omega_{J}$ and we obtain a second symplectic structure $\underline{\omega}_{J}$ on $\mathcal{C}/\mathfrak{g}_x$ (this can also be viewed as reduction in the Dirac setting \cite{BCWZ}
since $(S, \omega_{J})$ is a presymplectic realization \cite{BCWZ} of $M$ endowed with the Dirac structure (\ref{the-Dirac})). Then 
one considers $\underline{J}= \underline{\omega}^{-1}\underline{\omega}_{J}$, so that $J$ automatically commutes with 
$\underline{\omega}$ and $\underline{\omega}_{\underline{J}}$ is closed. Since $\omega+ J^*\omega$ is a pull-back by $\mu$, it is
zero when restricted to $\mathcal{C}$, and then we deduce that $\underline{J}^{2}= -Id$.
\end{proof}

\begin{remark}\rm \ Similar to the Poisson case (see e.g. \cite{CrFe02}), one can perform reduction starting with any generalized holomorphic map
$\mu: N\rmap M$ between two generalized complex manifolds. Then, for $x\in M$, we consider the isotropy Lie algebra $\mathfrak{g}_x$ at $x$, and we make the same assumptions as above ($x$ is a regular point and $\mathcal{C}/\mathfrak{g}_x$ is smooth). This time, the action of $\mathfrak{g}_x$ on $\mathcal{C}$ is $\alpha\mapsto \pi^{\sharp}\mu^*(\alpha)$, where $\pi$ is the bivector on $N$.
Then $\mathcal{C}/\mathfrak{g}_x$ will carry a natural Hitchin pair.
\end{remark}

\section{Bringing Lie algebroids/groupoids into the picture}

In this section we show that conditions (C1)-(C3) (see Proposition \ref{first}) have simple interpretations after passing to groupoids.
We start by recalling some basic facts on Lie groupoids and multiplicative forms.

\subsection{Lie algebroids/groupoids}
\label{Lie algebroids/groupoids:} Recall that a Lie algebroid over $M$ is a vector bundle $A\rmap M$ together with a Lie algebra bracket $[\cdot, \cdot]$ on the
space of sections $\Gamma(A)$ and a bundle map $\rho: A\rmap TM$, called the anchor of $A$, satisfying the Leibniz identity:
\[ [\alpha, f\beta]= f[\alpha, \beta]+ L_{\rho(\alpha)}(f)\beta ,\]
for all $\alpha, \beta\in \Gamma(A)$, $f\in C^{\infty}(M)$. In what follows, we will write $L_{\rho(\alpha)}$ simply as $L_{\alpha}$.
The basic examples are Lie algebras (when $M$ is a point) and tangent bundles (when $\rho$ is the identity). 

Lie groupoids are the global counterparts of Lie algebroids. In order to fix our notation, we recall
that a Lie groupoid over a manifold $M$ consists of a manifold $\grd$ together with
surjective submersions $\tar,\sour:\grd \to M$, called \textbf{target} and \textbf{source},
a partially defined multiplication $m: \grd^{(2)}\to \grd$, where
$\grd^{(2)}:=\{(g,h) \in \grd \times \grd\,|\, \sour(g)=\tar(h)\}$, a \textbf{unit section}
$\varepsilon:M \to \grd$ and an \textbf{inversion} $i:\grd \to \grd$, all related by the
appropriate axioms (see e.g. \cite{CrFe01} and the references therein). To simplify our notation, we will often
identify an element $x \in M$ with its image $\varepsilon(x) \in \grd$. We say that
$\grd$ is \textbf{source-simply-connected} if all the fibers $s^{-1}$ are connected and simply connected.

For a Lie groupoid $\grd$, the associated Lie algebroid $A(\grd)$
consists of the vector bundle
$$
\ker(d\sour)|_M \to M,
$$
with anchor
$\rho=d\tar: \ker(d\sour)|_M\to TM$ and bracket induced from the
Lie bracket on $\mathcal{X}(\grd)$ via the identification of
sections $\Gamma(\ker(d\sour)|_M)$ with right-invariant vector
fields on $\grd$ tangent to the $\sour$-fibers.

An \textbf{integration} of a Lie algebroid $A$ is a Lie groupoid
$\grd$ together with an isomorphism $A\cong A(\grd)$. Unlike Lie
algebras, not  every Lie algebroid admits an integration, see
\cite{CrFe01} for a description of the obstructions. On the other
hand, if a Lie algebroid is integrable, then there exists a
canonical source-simply-connected integration $\grd(A)$, and any other source-simply-connected integration
is smoothly isomorphic to $\grd(A)$;  see \cite{CrFe01}. To avoid endless repetitions, {\it from now on we assume 
that all Lie groupoids are source-simply-connected}.

Next, we recall a few facts about multiplicative forms. A form $\omega$ on a Lie groupoid $\grd$ is called
\textbf{multiplicative} if 
\begin{equation}\label{eq:compatible}
m^*\omega = \pr_1^*\omega + \pr_2^*\omega,
\end{equation}
where $\pr_i:\grd^{(2)}\to \grd$, $i=1,2$, are the canonical projections. 
Note that, for any form $\phi$ on $M$, $s^*\phi- t^*\phi$ is multiplicative;
we say that a multiplicative form \textbf{vanishes cohomologically} if it is of this type.
Cohomological vanishing is often as good as honest vanishing. According to this principle,
given a 3-form $\phi$ on $M$, one says that a 2-form $\omega$ on $\grd$ is $\phi$-closed if
$d\omega= s^*\phi- t^*\phi$.

On the infinitesimal side, for an arbitrary Lie algebroid $A$, one has the notion of \textbf{IM form on $A$}
(read infinitesimal multiplicative form): am IM form on $A$ is a bundle map
\[
u: A\rmap T^*M
\]
satisfying the following properties:
\begin{eqnarray}
\SP{u(\alpha),\rho(\beta)} & = & -\SP{u(\beta),\rho(\alpha)};\label{eq:skew}\\
u([\alpha, \beta]) & = & L_{\alpha}(u(\beta))-
                                      L_{\beta}(u(\alpha))+
                                      d\SP{u(\alpha),\rho(\beta')} ,\label{eq:Dirac}
\end{eqnarray}
for $\alpha,\beta \in \Gamma(A)$ (here $\SP{\cdot,\cdot}$ denotes the usual pairing between a 
vector space and its dual).

In general, if $A$ is the Lie algebroid of a Lie groupoid $\grd$, then a closed multiplicative 2-form $\omega$
on $\grd$ induces a IM form $u_{\omega}$ of $A$ by:
\[ \SP{u_{\omega}(\alpha), X}= \omega(\alpha, X) .\]
This has been explained in \cite{BCWZ} (Proposition 3.5), to which we refer for more details on multiplicative 2-forms. In particular, we 
will need the following reconstruction result (Theorem 5.1 in \cite{BCWZ}):

\begin{theorem}\label{Duke} 
If $A$ is an integrable Lie algebroid and if $\grd$ is its (source-simply-connected) integration, then 
$\omega\mapsto u_{\omega}$ is a 1-1 correspondence between closed multiplicative 2-forms on $\grd$ and 
IM forms of $A$.
\end{theorem}

A few more words on multiplicative structures on groupoids. Similar to 2-forms (and symplectic groupoids), one can talk
about the multiplicativity of various other structures. For instance, the multiplicativity of bivectors (and then Poisson groupoids) appears in \cite{MX}. Of relevance here are the $(1, 1)$-tensors: given a Lie groupoid $\grd$, we say that a $(1, 1)$-tensor $J: T\grd\rmap \grd$ is multiplicative if  for any 
$(g, h)\in \grd^{(2)}$ and any 
$v_g\in T_g\grd$, $w_h\in T_h\grd$ such that $(v_g, w_h)$ is tangent to $\grd^{(2)}$ at $(g, h)$, so is $(Jv_g, Jw_h)$,
and 
\[ (dm)_{g, h}(Jv_g, Jw_h)= J((dm)_{g, h}(v_g, w_h)).\]

With these, the various structures that we have discussed have a groupoid-version by requiring multiplicativity. For instance, 
a \textbf{symplectic groupoid} is a
Lie groupoid endowed with a symplectic form which is multiplicative.
Also, a \textbf{symplectic$+$complex groupoid}
is a Lie groupoid $\grd$ together with a symplectic$+$complex structure $(\omega, J)$ such that $\omega$ and $J$ are multiplicative.

\subsection{Taking care of (C1)} 

Important for us here is the relation of Poisson geometry with algebroids and symplectic groupoids
that we shortly recall. To a Poisson bivector $\pi$ on $M$, one associates a Lie algebroid structure on $T^*M$: 
the anchor is  $\pi^{\sharp}$, while the bracket is precisely the bracket $[\cdot, \cdot]_{\pi}$ defined in (\ref{br-pi}). 
Note that this bracket is
uniquely determined by the condition that $[df, dg]= d\{f, g\}$ and the Leibniz identity.

Now, (C1) says that $\pi^{\sharp}$ preserves the bracket, and this is easily seen to be equivalent to the fact that
$\pi$ is Poisson (actually, this is also equivalent 
to the fact that $[\cdot, \cdot]_{\pi}$ satisfies the Jacobi identity).
One says that $(M, \pi)$ is an \textbf{integrable Poisson manifold} (or that $\pi$ is integrable) if the induced algebroid is, and in this case
we denote by $\Sigma(M)= \Sigma(M, \pi)$ the associated source-simply connected Lie groupoid. Also, we say that a generalized complex structure is integrable
if the underlying Poisson manifold is. We will be interested in the integrable case only.

What is also important (but maybe too obvious to be noticed) is that the identity map
$Id_{T^*M}$ is an IM form on $T^*M$, hence, under the integrability assumption, it induces a multiplicative closed 2-form 
on $\Sigma(M)$ (cf. Theorem \ref{Duke}), that we denote by $\omega$. This is easily seen to be non-degenerate, hence 
$(\Sigma(M), \omega)$ is a symplectic groupoid. Conversely, if $\Sigma$ is a symplectic groupoid over $M$,
then there is an induced Poisson structure $\pi$ on $M$ uniquely determined by the condition that the source-map
is Poisson; moreover, the Lie algebroid  structure on $T^*M$ induced by $\pi$ is integrable (integrated by $\Sigma$). These constructions describe the well-known
correspondence (see \cite{CrFe02} and the references therein):

\begin{theorem} There is a natural 1-1 correspondence between:
\begin{enumerate}
\item[(i)] integrable Poisson structures on $M$ (i.e. $\pi$'s satisfying (C1), integrable).
\item[(ii)] symplectic groupoids $(\Sigma, \omega)$ over $M$.
\end{enumerate}
\end{theorem}

\subsection{Taking care of (C2)} 

\begin{theorem} 
\label{t-C2} Let $\pi$ be an integrable Poisson structure on $M$, and let $(\Sigma, \omega)$ 
be the associated symplectic groupoid over $M$. Then there is a natural 1-1 correspondence between
\begin{enumerate}
\item[(i)] $(1, 1)$-tensors $a$ on $M$ satisfying (C2).
\item[(ii)] multiplicative $(1, 1)$-tensors $J$ on $\Sigma$ with the property that $(J, \omega)$ is a Hitchin pair.
\end{enumerate}
In particular, there is a 1-1 correspondence between pairs $(\pi, a)$ on $M$ satisfying (C1) and (C2), with $\pi$ integrable,
and groupoids $\Sigma$ over $M$ together with a Hitchin pair $(\omega, J)$ on $\Sigma$ with $\omega$ and $J$ multiplicative.
\end{theorem}

\begin{proof}
With the bracket $[\cdot, \cdot ]_{\pi}$ and the IM conditions (equations (\ref{eq:skew}), (\ref{eq:Dirac})) in mind, the meaning of
(C2) is now clear: (C2) is equivalent to the fact that $a^{*}$ is an IM form on the algebroid $T^*M$ associated to the Poisson
structure $\pi$. Hence we obtain a 1-1 correspondence with closed multiplicative forms on $\Sigma$, i.e., cf.
Lemma \ref{lemma-complex+sympl} (see also the comment before Proposition \ref{nondeg}), with $(1, 1)$-tensors $J$ on $\Sigma$
with the property that $(\omega, J)$ is a twisted s.c. structure and $\omega_{J}$ is multiplicative. Finally, remark that
 Lemma \ref{lemma-complex+sympl} applied to symplectic groupoids $(\Sigma, \omega)$ has a bonus: $\omega_{J}$ is multiplicative if and only if
$J$ is.
\end{proof}

\subsection{Taking care of (C3)} 

We have seen that (C1) and (C2) alone
insure us that $\Sigma$ has an induced Hitchin pair $(\omega, J)$.
Note that the twist (i.e. the 2-form $-(\omega+ J^{*}\omega)$) is a multiplicative 2-form, but it may be non-zero. 
We will see below that what (C3) forces is that the twist vanishes cohomologically. This can be interpreted also using
generalized holomorphic maps. To motivate the statement recall that one of the main properties of symplectic groupoids $\Sigma$ over $M$ is
that the target map $t$ is Poisson while $s$ is anti-Poisson. Equivalently, $(t, s): \Sigma\rmap M\times \overline{M}$ is a Poisson map, where $\overline{M}$ is $M$ together with the opposite Poisson tensor $-\pi$. There are various other properties of groupoids
that can be expressed in terms of the map $(t, s)$. For instance, a groupoid $\Sigma$ is proper if $(t, s)$ is a proper map, or the presymplectic groupoids of \cite{BCWZ} have the main property that $(t, s)$ defines a Dirac realization of the base manifold.
In our context, if $M= (M, \mathcal{J})$ is a generalized complex manifold, then $\overline{M}$ is the opposite generalized complex manifold, i.e. $M$ together with $\overline{\mathcal{J}})$ (see (\ref{opp-J})).

\begin{theorem}
\label{t-C3}
Assume that $(\pi, a)$ satisfy (C1), (C2), with $\pi$ integrable, 
and let $(\Sigma, \omega, J)$ be the induced symplectic groupoid over $M$ together with the
induced multiplicative $(1, 1)$-tensor. Then, for a 2-form $\sigma$ on $M$, the following are equivalent:
\begin{enumerate}
\item[(i)] (C3) is satisfied.
\item[(ii)] $\omega+ J^{*}\omega = t^{*}\sigma- s^{*}\sigma$.
\item[(iii)] $(t, s): \Sigma \rmap M\times \overline{M}$ is generalized holomorphic.
\end{enumerate}
\end{theorem}

\begin{proof}
We first prove the equivalence between (i) and (ii).  
We consider the 2-form
\[ \phi= \tilde{\sigma}- t^{*}\sigma+ s^*\sigma,\]
where $\tilde{\sigma}= -\omega- J^{*}\omega$ is the twist of $(\omega, J)$ ($J= J_a$),
and let $A= Ker(ds)|_{M}$. 
We already know from Theorem \ref{Duke} that a closed multiplicative 2-form $\theta$ on $\Sigma$
is zero if and only if $u_{\theta}= 0$, i.e. 
\[ \theta(X_x, \alpha_x)= 0,\]
for all $X_x\in T_xM$, $\alpha_x\in A_x$. The proof of this fact is actually just a very simple step in
the proof of the theorem, and it appears as Corollary 3.4 in \cite{BCWZ}. We remark that exactly the same argument 
applies to higher degree forms. In particular, a 3-form $\theta$ is zero if and only if
\begin{equation}
\label{ad-hoc} 
\theta_x(X_x, Y_x, \alpha_x)= 0
\end{equation}
for all $X_x, Y_x\in T_xM$, $\alpha_x\in T_x\Sigma$. 
We can apply this to $d\phi$ (so that we can check when $\phi$ is closed) and then the original criteria 
to the 2-form $\phi$. Recalling also that $i_{X\wedge Y}(d\tilde{\sigma})= i_{\mathcal{N}_{J}(X, Y)}(\omega)$
(cf. the previous section on the non-degenerate case), we deduce that $\phi$ vanishes if and only if:
\begin{eqnarray}
 & \omega(\mathcal{N}_{J}(X_x, Y_x), \alpha_x)= (d\sigma)(X_x, Y_x, \rho(\alpha_x), \label{C3-1}\\
 & -\omega(X_x, \alpha_x)- \omega(J(X_x), J(\alpha_x))= \sigma(X_x, \rho(\alpha_x)),\label{C3-2}
\end{eqnarray}
for all $X_x, Y_x\in T_xM$, $\alpha_x\in A_x$. 
To make this more explicit, we need some remarks. First of all, after identifying $A$ with $T^*M$ we know 
(from the definition of $\omega$ and $\omega_{J}$) that:
\[ \omega(\alpha_x, X_x)= \alpha_x(X_x), \omega_{J}(\alpha_x, X_x)= \alpha(a(X)).\]
Note that the second equation shows that $J$ restricted to $T_xM$ equals to $a$, while restricted to $A_x= T_{x}^{*}M$ it equals to $a^{*}$. In particular, $J$ preserves the space of vector fields on $\Sigma$ tangent to $M$; since this space is closed 
under the Lie bracket, we deduce that $\mathcal{N}_{J}(X_x, Y_x)$ is tangent to $M$, and equals to $\mathcal{N}_{J}(X_x, Y_x)$,
for all $X_x, Y_x\in T_xM$. With these, we see that (\ref{C3-1}) and (\ref{C3-2}) say that, for all 
$X_x\in T_xM$ and $\alpha_x\in T_{x}^{*}M$,
\begin{eqnarray*}
 & - \alpha_x(\mathcal{N}_{a}(X_x, Y_x))= i_{X_x\wedge Y_x}(d\sigma)(\pi^{\sharp}(\alpha_x)) (= - \alpha_x(\pi^{\sharp}i_{X_x\wedge Y_x}(d\sigma)))\\
 & \alpha_x(X_x+ a^{2}X_x)= - \alpha_x(\pi^{\sharp}\sigma_{\sharp}(X_x)),
\end{eqnarray*}
and this is precisely (C3).

We now prove the equivalence of (ii) and (iii). Clearly, (ii) says that $(t, s)$ is compatible with the 2-forms. Also, $(t, s)$ is compatible with
the bivectors since $\Sigma$ is a symplectic groupoid. We claim that also the compatibility of $(s, t)$ with the $(1, 1)$-tensors is independent
of (C3) (i.e. it follows from the previous conditions (C1) and (C2)). Explicitly, this compatibility translates into the equations
\[ (dt)\circ J= a\circ (dt), \ \ (ds)\circ J= - a\circ (ds) .\]
We prove the first one (the second one being proven similarly, or deduced from the first one using $s= t\circ i$, $i(J(V))= - J(V)$). 
The main remark is that, for any multiplicative 2-form $\omega$ on any groupoid $\Sigma$, 
\[ \omega(\alpha, V)= \omega(\alpha, dt(V))= \langle u_{\omega}(\alpha), dt(V)\rangle,\]
for all $\alpha\in \Gamma(A)$ and all $V\in \mathcal{X}(\Sigma)$ (where $\alpha$ is extended to a vector field on $\Sigma$ using right translations).
This is one of the first consequences of the multiplicativity of $\omega$ (see equation (3.4) in \cite{BCWZ}). 
We apply this equation to our symplectic form $\omega$ (for which $u_{\omega}= Id$) and to the multiplicative form $\omega_{J}$ (for which $u_{\omega_{J}}= a^*$):
\[ \langle \alpha, a dt(V)\rangle = \langle a^{*}\alpha, dt(V)\rangle= \omega_{J}(\alpha, V)= \omega(\alpha, J(V))= \langle \alpha, dt(J(V))\rangle,\]
for all $\alpha \in T^*M$, hence the desired equation.
\end{proof}

\begin{definition} 
A \textbf{Hitchin groupoid} is a Lie groupoid $\Sigma$ together with a Hitchin pair $(\omega, J)$ on $\Sigma$ and a 2-form $\sigma$ on the the base such that $\omega+ J^{*}\omega= t^*\sigma- s^*\sigma$.
\end{definition}

With these, we conclude that there is a 1-1 correspondence between triples $(\pi, a, \sigma)$ satisfying (C1), (C2) and (C3), 
with $\pi$ integrable, and Hitchin groupoids $(\Sigma, \omega, J)$ over $M$. 
Note that an immediate consequence of the definition is that the isotropy groups $\Sigma_x= s^{-1}(x)\cap t^{-1}(x)$ are complex Lie groups with the complex structure coming from the restriction of $J$ (see also Corollary \ref{fibers}).

\subsection{A few more remarks} 

As we have seen throughout this paper, the Poisson and Dirac geometry underlying generalized complex structures can serve as a guide 
for various constructions and results in generalized complex geometry.   
We end with a few remarks along this line.

\begin{remark}\rm \
As it follows from the Poisson case, given a generalized complex manifold $M$, the associated Hitchin groupoid $\Sigma$ acts on any Hitchin realization $\mu: S\rmap M$ with $S$ compact (see e.g. \cite{CrFe02, BCWZ}). The action $\Sigma\times_{M}S\rmap S$ satisfies certain compatibility relations with the generalized complex structure (of course, one of them ensures that the action is Hamiltonian). 
In particular, the complex Lie group $\Sigma_x$ acts on the fiber $\mu^{-1}(x)$, and the reduced space $\mu^{-1}(x)/\Sigma_x$ is precisely the reduced space $\mathcal{C}/\mathfrak{g}_x$ discussed at the end of section \ref{Generalized holomorphic maps and reduction}.
\end{remark}

\begin{remark}\rm \
One can talk about dual pairs and Morita equivalence in the setting of generalized complex structures. Given two
generalized complex manifolds $M_1$ and $M_2$, a Morita equivalence will be given by a manifold $M$ together with a
Hitchin structure on it, and two submersions $\pi_i: M\rmap M_i$ ($i\in \{1, 2\}$) so that $M$ defines a Morita equivalence
between the underlying Poisson manifolds and $(\pi_1, \pi_2): M\rmap M_1\times\overline{M}_2$ is generalized holomorphic.
With these, the Hitchin groupoid of a generalized complex manifold $(M, \mathcal{J})$ defines a self Morita equivalence. 
\end{remark}

\begin{remark}\rm \ 
Recall \cite{Gu} that, given a generalized complex structures $\mathcal{J}$ on $M$, and a closed 2-form $B$, the gauge transformation
of $\mathcal{J}$ with respect to $B$ is the new generalized complex structure  $\mathcal{J}_{B}$ given by the matrix
\begin{equation}
\mathcal{J}_B= \left( \begin{array}{ll}
                1 & 0 \\
                -B & 0  
        \end{array}
   \right)  
\mathcal{J}
\left( \begin{array}{ll}
                1 & 0 \\
                B & 0  
        \end{array}
   \right) 
\end{equation}
One says that $\mathcal{J}$ and $\mathcal{J}'$ are gauge equivalent if $\mathcal{J}'= \mathcal{J}_B$ for some closed 2-form $B$. 
As in \cite{BO} (see Theorem 5.1 there), if $(M, \mathcal{J})$ is integrable and $(\Sigma, \omega, J)$ is the associated Hitchin groupoid,
then one can recover the Hitchin groupoid of $(M, \mathcal{J}_{B})$ and conclude that $(M, \mathcal{J})$ and $(M, \mathcal{J}_{B})$
are Morita equivalent. More precisely, one defines the gauge transformation of $J$ as the new $(1, 1)$-tensor on $\Sigma$ given by
\[ J_{B}= J+ \omega^{-1}(s^*B-t^*B) ,\]
and then $(\Sigma, \omega, J_{B})$ is the Hitchin groupoid of $(M, \mathcal{J}_{B})$. Moreover, $(\Sigma, \omega, J- \omega^{-1}t^{*}(B))$
defines a Morita equivalence between $(M, \mathcal{J})$ and $(M, \mathcal{J}_B)$.
\end{remark}

\begin{remark} Some of the remarks above also follow from the corresponding results for Dirac structures and the fact that, if
$(\Sigma, \omega, J)$ is the Hitchin groupoid of a generalized complex structure $\mathcal{J}$, then $(\Sigma, \omega_{J})$ is the presymplectic groupoid \cite{BCWZ} of the Dirac structure (\ref{the-Dirac}) underlying $\mathcal{J}$.
\end{remark}

\begin{footnotesize}

\end{footnotesize}

\end{document}